\documentclass[11pt]{article}
\usepackage{geometry}                
\usepackage[parfill]{parskip}    
\usepackage{graphicx}
\usepackage{amsmath, latexsym, amssymb}
\usepackage{epstopdf}
\input xypic
\begin{document}
\title{ Universal spaces for    manifolds equipped with a closed  integral k-form }
\author{\bf  H\^ong-V\^an L\^e \footnote{Mathematical Institute  of ASCR,
Zitna 25, Praha 1, CZ-11567 Czech Republic,
  email: hvle@math.cas.cz}}
\date{}
\maketitle
\newcommand{\R}{{\mathbb R}}
\newcommand{\C}{{\mathbb C}}
\newcommand{\F}{{\mathbb F}}
\newcommand{\Z}{{\mathbb Z}}
\newcommand{\N}{{\mathbb N}}
\newcommand{\Q}{{\mathbb Q}}
\newcommand{\Hq}{{\mathbb H}}

\newcommand{\Aa}{{\mathcal A}}
\newcommand{\Bb}{{\mathcal B}}
\newcommand{\Cc}{{\mathcal C}}    
\newcommand{\Dd}{{\mathcal D}}
\newcommand{\Ee}{{\mathcal E}}
\newcommand{\Ff}{{\mathcal F}}
\newcommand{\Gg}{{\mathcal G}}    
\newcommand{\Hh}{{\mathcal H}}
\newcommand{\Kk}{{\mathcal K}}
\newcommand{\Ii}{{\mathcal I}}
\newcommand{\Jj}{{\mathcal J}}
\newcommand{\Ll}{{\mathcal L}}    
\newcommand{\Mm}{{\mathcal M}}    
\newcommand{\Nn}{{\mathcal N}}
\newcommand{\Oo}{{\mathcal O}}
\newcommand{\Pp}{{\mathcal P}}
\newcommand{\Qq}{{\mathcal Q}}
\newcommand{\Rr}{{\mathcal R}}
\newcommand{\Ss}{{\mathcal S}}
\newcommand{\Tt}{{\mathcal T}}
\newcommand{\Uu}{{\mathcal U}}
\newcommand{\Vv}{{\mathcal V}}
\newcommand{\Ww}{{\mathcal W}}
\newcommand{\Xx}{{\mathcal X}}
\newcommand{\Yy}{{\mathcal Y}}
\newcommand{\Zz}{{\mathcal Z}}

\newcommand{\zt}{{\tilde z}}
\newcommand{\xt}{{\tilde x}}
\newcommand{\Ht}{\widetilde{H}}
\newcommand{\ut}{{\tilde u}}
\newcommand{\Mt}{{\widetilde M}}
\newcommand{\Llt}{{\widetilde{\mathcal L}}}
\newcommand{\yt}{{\tilde y}}
\newcommand{\vt}{{\tilde v}}
\newcommand{\Ppt}{{\widetilde{\mathcal P}}}

\newcommand{\Remark}{{\it Remark}}
\newcommand{\Proof}{{\it Proof}}
\newcommand{\ad}{{\rm ad}}
\newcommand{\Om}{{\Omega}}
\newcommand{\om}{{\omega}}
\newcommand{\eps}{{\varepsilon}}
\newcommand{\Di}{{\rm Diff}}
\newcommand{\Pro}[1]{\noindent {\bf Proposition #1}}
\newcommand{\Thm}[1]{\noindent {\bf Theorem #1}}
\newcommand{\Lem}[1]{\noindent {\bf Lemma #1 }}
\newcommand{\An}[1]{\noindent {\bf Anmerkung #1}}
\newcommand{\Kor}[1]{\noindent {\bf Korollar #1}}
\newcommand{\Satz}[1]{\noindent {\bf Satz #1}}

\newcommand{\gl}{{\frak gl}}
\renewcommand{\o}{{\frak o}}
\newcommand{\so}{{\frak so}}
\renewcommand{\u}{{\frak u}}
\newcommand{\su}{{\frak su}}
\newcommand{\ssl}{{\frak sl}}
\newcommand{\ssp}{{\frak sp}}

\newcommand{\Cinf}{C^{\infty}}
\newcommand{\CS}{{\mathcal{CS}}}
\newcommand{\YM}{{\mathcal{YM}}}
\newcommand{\Jreg}{{\mathcal J}_{\rm reg}}
\newcommand{\Hreg}{{\mathcal H}_{\rm reg}}
\newcommand{\SP}{{\rm SP}}
\newcommand{\im}{{\rm im}}

\newcommand{\inner}[2]{\langle #1, #2\rangle}    
\newcommand{\Inner}[2]{#1\cdot#2}
\def\NABLA#1{{\mathop{\nabla\kern-.5ex\lower1ex\hbox{$#1$}}}}
\def\Nabla#1{\nabla\kern-.5ex{}_#1}

\newcommand{\half}{\scriptstyle\frac{1}{2}}
\newcommand{\p}{{\partial}}
\newcommand{\notsub}{\not\subset}
\newcommand{\iI}{{I}}               
\newcommand{\bI}{{\partial I}}      
\newcommand{\LRA}{\Longrightarrow}
\newcommand{\LLA}{\Longleftarrow}
\newcommand{\lra}{\longrightarrow}
\newcommand{\LLR}{\Longleftrightarrow}
\newcommand{\lla}{\longleftarrow}
\newcommand{\INTO}{\hookrightarrow}

\newcommand{\Sy}{\text{ Diff }_{\om}}
\newcommand{\Ex}{\text{Diff }_{ex}}
\newcommand{\jdef}[1]{{\bf #1}}
\newcommand{\QED}{\hfill$\Box$\medskip}

\newcommand{\UuU}{\Upsilon _{\delta}(H_0) \times \Uu _{\delta} (J_0)}
\newcommand{\bm}{\boldmath}

\medskip

\begin{abstract}
In this note we  prove that  any integral closed k-form $\phi ^k$, $k\ge 3$,
on a m-dimensional manifold $M^m$, $m \ge k$, is  the restriction of a  universal  closed  k-form $h^k$
on a universal manifold $U^{d(m,k)}$ as a result of an  embedding  of $M^m$
to $U^{d(m,k)}$.

\end{abstract}

\medskip


\section {Introduction.}

Let a manifold $M^m$  be equipped with a  tensor of degree $\alpha^k$ and a manifold $N^n$ 
 equipped with a tensor 
 $\beta^k$.  Suppose that $n> m$. We want to know if there is an immersion
$f: M^m \to N^n$ such that $f^*(\beta^k) = \alpha^k$. This problem has  a long history,  the  Nash  embedding theorem  for $(M^m,\alpha^2)$  being a  Riemannian   manifold and $(N^n, \beta^2)$ being  the standard  Euclidean
space is one of most spectacular  results in this field.
Gromov in his seminal book [Gromov1986]  developed many methods for solving  this problem. 

In this note we apply  the Gromov theory to  obtain the existence of a universal space
$U ^{d(m,k)}$,  equipped with an integral closed  k-form $h^k$, $k\ge 3$,  for any m-dimensional
manifolds $M^m$ equipped with an integral  closed k-form $\phi^k$ (Theorem 3.6).   Theorem 3.6
is a generalization of Tischler's theorem [Tischler1977] on the existence of a symplectic embedding  from an integral  symplectic manifold
$(M^{2n}, \om )$ to $\C P^n$ equipped with  the standard Kaehler symplectic form.
\medskip

This note also contains an Appendix written in communication with Kaoru Ono which
contains a new ``soft" proof
of a version Theorem 3.6 on the existence of a universal space for  manifolds equipped with an integral closed k-form.  Our soft proof
does not use the Nash-Gromov implicit function theorem, but  we do not get  a $C^0$-perturbation 
result as in Theorem 3.6.

\medskip
\section {  H-principle and Nash-Gromov  implicit
function theorem.}

In this section we briefly  recall some important notions and results in the Gromov theory [Gromov1986]
which we shall use for our proof of Theorem 3.6.

Let $V$ and $W$ be smooth manifolds. We denote by $(V,W)^{(r)}$, $r\ge 0$, 
 the space of $r$-jets  of smooth mappings from $V$ to $W$.
We shall think of each map $f : V \to W$ as a section of the fibration $V \times W
=(V, W)^{(0)}$
over $V$. Thus $(V, W) ^{(r)}$ is a fibration over $V$, and we shall denote
by $p ^r $ the canonical projection $(V, W)^{(r)}$ to $V$, and by
$ p^s_r$ the canonical projection $(V, W) ^{(s)} \to (V, W) ^{(r)}$, for any $s > r$.

A section $s : V \to (V,W)^r$ is called {\bf holonomic}, if  $s$ is the r-jet of some
section $f: V \to (V,W)$.

We   say that a  differential relation $\Rr\subset (V,W)^{(r)}$ satisfies the {\bf H-principle}, if  every continuous section $\phi_0: V \to \Rr$ can be brought to a holonomic section $\phi_1$
by a homotopy of sections  $\phi_t: V \to \Rr$, $t\in [0,1]$.

We   say that a  differential relation $\Rr\subset (V,W)^{(r)}$ satisfies the {\bf H-principle
$C^0$-near a map} $f_0 :  V \to W$, if  every continuous section $\phi_0: V \to \Rr$ which lies over $f_0$,
(i.e. $p ^r_0 \circ \phi_0 = f_0$) can be brought to a holonomic section $\phi_1$
by a homotopy of sections  $\phi_t: V \to \Rr _U$, $t\in [0,1]$, for an arbitrary
small neighborhood $U$ of $f_0(V)$ in $V\times W$ [Gromov1986, 1.2.2].  Here for an
open set 
$U\subset V\times W $, we write
$$ \Rr_U := (p ^r_0) ^{-1} (U) \cap \Rr \subset (V, W) ^r.$$

 The H-principle is called
{\bf $C^0$-dense}, if it holds true $C^0$-near every map $f: V \to W$.

\medskip

We also define  the {\bf   fine $C^0$-topology}  on the space
$C^0(X)$   of continuous  sections  of a smooth fibration $X \to V$ by taking  the sets
$C^0(U) \subset C^0 (X)$, $U$ is  open in $X$,  as the basis for 
this topology.  The {\bf fine $C^r$-topology on $ C^r (X)$} is induced by  the fine $C^0$-topology 
on $C^0(X^{(r)})$ using the embedding $C^r (X)\to C^0 ( X^{(r)})$.

\medskip

Suppose  we are given  a differential relation $\Rr\subset (V, W) ^{(r)}$. 
We define  the {\bf prolongation} $\Rr^k  \subset (V,W)^{r+k}$ inductively.  Let $\Rr' \subset (\Rr ^{(r)}) ^{(1)}$ consist of the 1-jets of germs  of $C^1$-sections $V\to \Rr$.
We put $\Rr^1: =  \Rr' \cap X^{(r+1)}\subset (X^{(r)})^{(1)}$.  Then repeat this and define
$\Rr^k :=  (V,W)^{r+k}\cap (\Rr^{k-1})^1\subset ((V,W)^{r+k-1})^1$. 
A $C^{r+k}$-solution
of $\Rr$ is a holonomic section of $\Rr^k$.

Fix an integer $k\ge r$ and denote
by $\Phi (U)$  the space of $C^k$-solutions of $\Rr$ over $U$ for all open $U\subset V$.
This set  equipped with the natural restriction $\Phi (U) \to \Phi (U')$ for
all $U' \subset U$ makes $\Phi$ a sheaf  which we call the {\bf solution sheaf} of $\Rr$ over $V$.
We shall say that $\Phi$ satisfies the
$H$-principle, if $\Rr$ satisfies the H-principle.

A sheaf $\Phi$ is called {\bf flexible (microflexible)}, if the restriction
map $\Phi (C) \to \Phi (C')$ is a fibration (microfibration) for all pair
of compact subsets $C$ and $C'\subset C$ in $M$. We recall that the map
$\alpha : A \to A'$ is called {\bf microfibration}, if the  homotopy  lifting
property  for a homotopy $\psi: P \times [0,1]\to A'$ is  valid only ``micro", i.e.
there exists $\eps > 0$ such that $\psi$ can lift to a homotopy  $\bar \psi : P\times [0,\eps]
\to A$.
\medskip

 \textbf{2.1. H-principle and flexibility} [Gromov1986, 2.2.1.B]. {\it If $V$ is a locally compact countable polyhedron (e.g. manifold), then every flexible sheaf over $V$ satisfies  the H-principle.}

\medskip

{\bf 2.2. A  criterion on flexibility.} [Gromov1986, 2.2.3.C''] {\it 
Let $\Phi$ be a microflexible sheaf over $V$ and let a submanifold $V_0 \subset V$ be
sharply movable by acting  diffeotopies.  Then the sheaf $\Phi_0 = \Phi_{|V_0}$ is flexible and hence it satisfies the H-principle.}

\medskip

One of Gromov's method to get the microflexibility of some sheaf  (and then to get the H-principle)  is to exploit the Nash-Gromov implicit function theorem.

 Let $X\to V$ be a smooth fibration and  $G\to V$ be a smooth vector bundle over a manifold $V$. We denote by $\Xx ^\alpha$ and $\Gg ^\alpha$ respectively the spaces  of $C^\alpha$-sections of the fibrations $X$ and 
$G$ for all $\alpha = 0,1,\cdots , \infty$.  Let $\Dd: \Xx^r\to \Gg^0 $ be a differential operator of order $r$. In other words the operator $\Dd$ is given by a  bundle map $\triangle : X^{(r)} \to G$, namely $\Dd(x) = \triangle \circ J^r_x$, where
$J^r _x(v)$ denotes the r-jet of $x$ at $v\in V$. We assume below that $\Dd$ is 
a $C^\infty$-operator and so we have continuous maps $\Dd : \Xx^{\alpha +r} \to \Gg^\alpha$ for all $\alpha = 0, 1, \cdots, \infty$.

 Now we shall define the linearization of  a differential operator $\Dd$. Let $x$ be a $C^{\alpha}$-section 
 of a smooth  vector bundle $X \to V$. Denote by $Y_x$ the induced  vector bundle $x^* (T_{vert} (X))$.
 For each $\beta \le \alpha$ we denote by  $\Yy_x ^\beta$ the space of $C^\beta$-section $V \to Y_x$.
 The space $\Yy_x ^\alpha$ can be  considered as the tangent space $T_x (\Xx^\alpha)$.
 Now  we suppose that the   fibration $X\to V$  does not have  boundary.  For  $x\in \Xx^r$ the
 {\bf linearization} $L_x : \Yy _x ^r \to \Gg ^0$  of the operators $\Dd$ at $x$ is defined as follows.
 Let $y = \p x_t / \p t_{|t = 0}$. Then
 $$ L_x (y) = L(x,y)  = {\p \over \p t} _{| t = 0} \Dd (x_t).$$
 We say that the operator $\Dd$  is {\bf infinitesimal invertible over a subset
$\Aa$} in the space of sections $x : V \to X$ if there exists a family of linear
differential operators of certain order $s$, namely $M_x: \Gg^s \to \Yy ^0_x$, 
for $x\in \Aa$, such that the following three properties are satisfied.

\begin{enumerate}
	\item There is an integer $d\ge r$, called {\bf the defect of the infinitesimal inversion $M$}, such that $\Aa$ is contained in $\Xx^d$, and  furthermore, $\Aa = \Aa^d$ consists exactly  of $C^d$-solutions of an open differential relation
	$A \subset X^{(d)}$. In particular, the sets $\Aa ^{\alpha +d} = \Aa \cap \Xx ^{\alpha +d}$ are open in $\Xx^{\alpha +d}$ in the respective fine $C^{\alpha +d}$-topology for all $\alpha = 0,1, \cdots, \infty$.
	\item The operator $M_x (g) = M (x, g)$ is a (non-linear) differential operator
	in $x$ of order $s$. Moreover the global operator
	$$ M : \Aa^d \times \Gg^s \to \Jj ^0 = T (\Xx^0)$$
	is a differential  operator, that is given by a $C^\infty$-map $A \oplus G^{(s)} \to T_{vert} (X)$.
	\item $L_x \circ M_x = Id$ that is 
	$$ L(x, M(x,g)) = g \text { for all } x\in \Aa ^{d+r} \text { and } g\in \Gg ^{r+s}.$$
\end{enumerate}

Now let $\Dd$ admit over an open set $\Aa = \Aa^d \subset \Xx^d$ an infinitesimal inversion $M$ of order $s$ and of defect $d$. 
For a subset $\Bb \subset \Xx^0 \times \Gg^0$ we put
$\Bb^{\alpha,\beta} : = \Bb \cap (\Xx^\alpha \times \Gg^\beta)$. 
Let us fix an integer $\sigma_0$ which satisfies the following inequality
$$\sigma_0 > \bar s = \max (d, 2r +s).\leqno (*)$$

Finally we fix an arbitrary Riemannian metric in the underlying manifold $V$.

\medskip

{\bf 2.3. Nash-Gromov implicit function theorem.} [Gromov1986, 2.3.2]. {\it  There exists a family of sets $\Bb _x\subset \Gg ^{\sigma_0 +s}$ for all $x
\in \Aa^{\sigma_0+r+s}$, and a family of operators $\Dd _x ^{-1} : \Bb_x \to \Aa$ with the following five properties.
\begin{enumerate}
	\item Neighborhood property: Each  set $\Bb _x$ contains a neighborhood of zero
	in the space  $\Gg ^{\sigma_0 +s}$. Furthermore, the union  $\Bb = \{x\} \times \Bb_x$ where $x$ runs over $\Aa ^{\sigma_0 +r+s}$, is an open subset in the space $\Aa ^{\sigma _0 +r +s}\times \Gg^{\sigma_0 +s}$.
	\item Normalization Property: $\Dd _x ^{-1} (0) = x$ for all $x\in \Aa ^{\sigma_0 +r+s}$.
	\item Inversion Property: $\Dd \circ \Dd_x ^{-1} - \Dd (x) = Id$, for all $x\in \Aa ^{\sigma _0+r+s}$, that is
	$$ \Dd ( \Dd _x ^{-1} (g)) = \Dd (x) + g,$$
	for all pairs $(x,g) \in \Bb$.
	\item Regularity and Continuity: If the section $x\in \Aa$ is $C^{\eta_1 +r+s}$-smooth and
	if $g\in \Bb_x$ is $C^{\sigma_1 +s}$-smooth for $\sigma_0 \le \sigma_1 \le \eta_1$, then the section $\Dd _x ^{-1} (g)$ is $C^\sigma$-smooth for all $\sigma < \sigma_1$. Moreover the operator $\Dd ^{-1} : \Bb^{\eta _1 + r +s, \sigma_1 +s} \to\Aa^\sigma$,
	$\Dd ^{-1} (x,g) = \Dd^{-1} _x (g)$, is jointly continuous  in the variables $x$ and $g$. Furthermore, for $\eta _1  > \sigma _1$, the section $\Dd^{-1} : \Bb ^{\eta_1 +r+s, \sigma_1 +s} \to  \Aa^{\sigma_1}$ is continuous.
	\item  Locality: The value of the section  $\Dd ^{-1} _x (g) : V \to X$ at any given
	point $v\in V$ does not depend on the behavior of $x$ and $g$ outside the unit ball
	$B_v (1)$ in $V$ with center $v$, and so the equality $(x,g) _{| B_v (1)} = (x', g')_{| B_v(1)} $ implies  $\Dd _x ^{-1} (g)) (v) = ( \Dd_{x'} ^{-1} (g') ) (v)$.
\end{enumerate}
}

\medskip

{\bf 2.4. Corollary. Implicit Funtion Theorem.} {\it For every $x_0\in \Aa ^\infty$ there exists  fine $C^{\bar s +s+1}$-neighborhood $\Bb_0$ of zero in the space of $\Gg_{\bar s+s +1}$, where $\bar s = \max (d, 2r+s)$, such that for each $C^{\sigma + s}$-section $g\in \Bb_0$, $\sigma \ge \bar s +1$, the equation $\Dd(x) = \Dd(x_0) +g$
has a $C^\sigma$-solution.}

\medskip

Finally we shall show a large class of  microflexible solution sheafs $\Phi$ by using
 the Nash-Gromov implicit function theorem.

 Let us fix a $C^\infty$-section $g: V \to G$ and we call  a $C^\infty$-germ $x: \Oo p (v) \to X$, $v\in V$, {\bf an infinitesimal solution of order $\alpha$ of the equation
 $\Dd (x) = g$}, if at the point $v$ the germ $g' = g - \Dd (x)$ has  zero $\alpha$-jet
 , i.e. $J ^\alpha_{g'} (v) = 0$.  We denote by $\Rr^\alpha (\Dd, g)\subset X^{(r+\alpha)}$ the set of all jets  represented by these infinitesimal  solutions of
 order $\alpha$ over all points $v\in V$.  Now we recall the open set $A\subset  X^{(d)}$ defining the set $\Aa \subset X^{(d)}$, and for $\alpha \ge d-r$ we put
 $$\Rr _\alpha = \Rr _\alpha (A, \Dd, g) = A ^{r+\alpha -d} \cap \Rr^\alpha (\Dd, g) \subset X^{(r+\alpha)},$$
 where $A^{r+\alpha -d} = ( p ^{r+\alpha} _d ) ^{-1} (A)$ for $p ^{r+\alpha }_d : X^{r+\alpha} \to X^d$.
 
A $C^{r+\alpha}$-section $x : V \to X$ satisfies $\Rr_\alpha$, iff $\Dd(x) = g$ and $x\in \Aa$.

Now we set $\Rr = \Rr_{d-r}$ and denote by $\Phi = \Phi (\Rr) = \Phi (A, \Dd, g)$ the sheaf of $C^\infty$-solutions of  $\Rr$.

\medskip

{\bf 2.5. Microflexibility of the sheaf of solutions and the Nash-Gromov implicit functions.}[Gromov1986 2.3.2.D''] {\it  The sheaf $\Phi$ is microflexible.  }

\medskip

\section { Universal space for integral closed k-forms on m-dimensional manifolds.}

Suppose that $m \ge  k \ge 3$. In this section we shall show that any integral closed k-form $\phi^k$ on a m-dimensional
smooth manifold $M^m$ can be induced from a  universal  closed k-form $h^k$  on a universal manifold $U^{d(m,k)}$  by an embedding $M^m$ to   a universal space $( U^{d(m,k)},  h^k)$, see Theorem 3.6.

\medskip

Our definition of the universal space $( U^{d(m,k)},  h^k)$ is based on the work of
Dold and Thom [D-T1958] as well as the idea of Gromov [Gromov2006] to reduce this
problem to the case  that $\phi$ is an exact $k$-form.

\medskip

Let $SP^q(X)$ be the $q$-fold symmetric product of a locally compact, paracompact
Hausdorff pointed space $(X,0)$ , i.e. $SP^q(X)$ is the quotient space of the $q$-fold Cartesian
$(X^q,0)$ over the  permutation group $\sigma _q$. We shall denote
by $SP (X,0)$ the  inductive limit of $SP^q(X)$ with the inclusion
$$X =SP^1(X)  \stackrel{i_1}{\to}SP^2(X) \stackrel{i_2}{\to} \cdots \to 
SP ^q(X)\stackrel{i_q}{\to} \cdots,  $$
where
$$ SP^q (X) \stackrel{i_q}{\to} SP^{q+1}(X): \: [x_1,x_2, \cdots , x_q]\mapsto [0, x_1,x_2, \dots , x_q].$$

Equivalently we can write 
$$SP (X,0) = \sum_q SP^q (X) / ([x_1, x_2, \cdots , x_q]
\sim [0, x_1, x_2, \cdots, x_q]).$$

So we shall also denote by $i_q$ the canonical inclusion $SP ^q(X) \to SP (X,0)$.

\medskip

{\bf 3.1. Theorem.} (see [D-T1958, Satz 6.10])  {\it There exist natural isomorphisms
$j : H_q(X, \Z) \to \pi_q (SP(X,0))$   for $q> 0$.}

\medskip

{\bf 3.2. Corollary.} ([D-T1958]) {\it The space $SP (S^n,0)$ is the Eilenberg-McLane
complex $K(\Z, n)$.}

\medskip

Now let $\tau^k$  be the generator of $H^k (SP (S^k, 0),\Z)$ and by abusing notations we also denote by $\tau^k$ the restriction of the generator $\tau^k$ to any subspace $i_q(SP^q (S^k)) \subset SP (S^k, 0)$. The following lemma shows that we can replace a classifying map from $(M^m,[\phi^k])$ to
$(SP (S^k, 0), \tau ^k)$ by a map from $(M^m, [\phi^k])$ to $(SP ^{ [{m-k\over 2}] +1} (S^k), \tau ^k)$ .

\medskip

{\bf 3.3. Lemma.} {\it Let $[\phi^k] \in H^k (M^m, \Z)$. Then there exists a continuous map  $ f $ from $M^m$ to  $(SP^{ [{m-k\over 2}] +1}(S^k))$ such that $f^* (\tau ^k) = [\phi^k]$.}

\medskip

{\it Proof.} Let $f_0$ be a classifying map  from $M^m$ to $SP (S^k, 0)$ such that
$f_0 ^* (\tau^k) = \alpha$. Denote by $K^i$ the i-dimensional skeleton of $SP (S^k, 0)$ and by
$\bar \tau^k$ the restriction  of $\tau^k$ to $K^m$. Then we know that $f_0$ is homotopic equivalent  to a continuous map $f_1 : M^m \to K^m$ such that $f_1 ^* (\bar \tau^k) = [\phi^k]$. To prove Lemma 3.3 it suffices to find a map $g: K^m \to SP^{ [{m-k\over 2}] +1} (S^k)$ such that $g^* (\tau^k) = \bar \tau^k$. Then the map $f = g\circ f_1$ satisfies the  condition of Lemma 3.3.

We observe that $K^{k+1} = K^k$ consists of  the sphere $S^k$. If $m= k$ or $m = k +1$, then $g$ can be chosen as  the identity map.  Now suppose that $m \ge  k +2$.  The  following identity [D-P1961, (12.12)] 
 
$$\pi _i (SP^n (X)) = H_i (X) \text { for } i < k +2n -1, \: n > 1,$$
if $X$ is connected  and $H_i (X) = 0$ for $ 0< i < k, \, k > 1,$
implies  that 
 $$\pi_i( SP^{ [{m-k\over 2}] +1} (S^k))= 0, \text { for }  k+1\le  i \le m-1. $$

Using the obstruction theory we obtain a map $g : K^m \to SP^ {[{m-k\over 2}] +1}(S^k)$ extending the inclusion  map
$K^k = S^k\to  SP^ {[{m-k\over 2}] +1}(S^k)$. Clearly the map $g$ satisfies the
required property that $g^* (\tau^k) = \bar \tau^k$.
\QED

\medskip

Since $SP^{[{m-k\over 2}] +1}(S^k)$ has a finite simplicial decomposition we can apply the Thom construction in [Thom1954, III.2] where Thom showed that  any finite m-dimensional polyhedron $K$ can be
embedded in  a  compact $(2m+1)$-dimensional manifold $M^{2m+1}$ such that $K$ is
a retract of $M^{2m+1}$.  As a result we get the following
\medskip

{\bf 3.4. Lemma.}  {\it The   space $SP^{[{m-k\over 2}] +1}(S^k)$  can be embedded into a compact  smooth manifold $\Mm^{s(m,k)}$, $s(m,k) = 2 ([{m-k\over 2}] +3)k +1$,
such that  (the image of) $SP^{[{m-k\over 2}] +1}(S^k)$ is a  retract of  $\Mm^{s(m,k)}$.}

\medskip

Let us denote  also by $\tau^k$ the pull back of the universal class $\tau^k$ from $SP^{[{m-k\over 2}] +1} (S^k)$ to $\Mm^{s(m,k)}$ and let $\alpha^k$ be any differential form representing $\tau^k$
on $\Mm^{s(m,k)}$.

Let $\beta^k_l$ be the following k-form on $\R^{k\cdot l}$ with coordinates $x^{ij}, 1\le i \le l, \, 1\le j \le k$
$$\beta_l^k = dx^{11} \wedge dx^{12} \cdots \wedge dx^{1k} + \cdots  +dx^{l1} \wedge dx^{l2} \cdots \wedge dx^{lk}.$$

\medskip

Set $d(m,k) := s(m,k) + 2m +2 - k + {1\over 2} (k-1) [ {2m\over k-1}]([ {2m\over k-1}]-1) +
k (m+1) \binom {m+1} {k} $.

\medskip

Now we state the main theorem of this section. 
Let
$$(U^{d(m,k)},h^k) = (\Mm^{s(m,k)}\times \R^{kN}, \alpha^k \oplus \beta^k_{N}). \leqno(3.5)$$

\medskip

{\bf 3.6. Theorem.} {\it Suppose that $\phi^k$ is a closed integral k-form
on a  smooth manifold $M^m$. Then there exists an embedding $f : M^m\to (U^{d(m,k)}, h^k)$ such that $f^* (h^k) = \phi$. Moreover for any given map $\tilde f : M^m \to  ( U^{d(m,k)},  h^k)$ such that $\tilde f^* [h^k] = [\phi^k]$ there exists a $C^0$-close to $\tilde f$ embedding $f : M^m \to U^{d(m,k)}$ 
such that $f ^* ( h^k) = \phi^k$.}

\medskip

{\it Proof of Theorem 3.6}. Using Lemma 3.3 and Lemma 3.4 we see that  the first statement of Theorem 3.6 follows from the second statement of Theorem 3.6. Furthermore we shall reduce  the second statement to an immersion problem for
exact 3-forms as follows.  Denote by $\tilde f_1: M^m \to \Mm^{s(m,k)}$  the projection of $\tilde f$ to the first factor. Then we have
$\tilde f_1 ^* (\tau^k) = [\phi^k] \in H^k (M^m,\Z)$. Let
$$ g= \phi - \tilde f_1^* (\alpha).$$
Clearly $g$ is an exact k-form on $M^m$. We can also assume that $\tilde f_1$ is an embedding
by  perturbing this map a little,  if necessary.  Thus the second statement of Theorem 3.6 is a corollary of
the following Proposition.

\medskip

{\bf  3.7. Proposition.} {\it For any given map $f_0 : M^m \to \R^{kN}$ there is a
$C^0$-close to $f_0$ immersion $f_3 : M^m \to  (\R^{kN}, \beta^k_{N})$ such that $f_3 ^* (\beta^k_{N}) = g$ for any  exact  $k$-form $g$.}

\medskip

 We shall apply  the Gromov H-principle for immersion of  differential forms to prove Proposition 3.7. Gromov  extended the Nash idea to add  some regularity  for an immersion in order
to  apply the implicit function  theorem and  then  using 2.5 to get the H-principle for  the isometric immersion.  Finally   using the H-principle we shall  get  immersions required
in Proposition 3.7.

 Let
$h$ be a smooth differential $k$-form on a manifold $W$.  Denote by $I_{h(w)}$ a linear homomorphism
$$ I_{h(w)} : T_wW \to  \Lambda ^{k-1} (T_wW)^*, \:   X \mapsto  X \rfloor  h.$$
A subspace  $T \subset T_w W$ is called {\bf $h(w)$-regular}, if the composition
of $I_{h(w)}$ with the restriction homomorphism  $r: \Lambda ^{k-1} (T_wW)^* \to \Lambda ^{k-1}( T)^*$
sends $T_wW$ onto $\Lambda ^{k-1}( T)^*$.

An immersion $f: V \to W$  is called {\bf h-regular},
if  for all $v \in V$  the subspace $Df (T_v V)$ is $h(f(v))$-regular.
 
\medskip

{\it Proof of Proposition 3.7.} Roughly speaking, we 
add the  condition of  $\beta^k_{N}$-regularity to the isometry property (i.e. $f_3^* ( \beta^k_{N}) = g$) and
extend this equation for  mappings also denoted by $f_3$ from the manifold
$M ^{m+1} = M^m \times (-1, 1)$ provided with a form $ g \oplus 0$, denoted from  now on also  by $g$, to  the space $(\R^{kN}, \beta^k_{N})$.  Our Proposition 3.10 states
that  the solution sheaf restricted to $M^m\subset M^{m+1}$ satisfies the H-principle. In fact, this statement is
a consequence of  Theorem 3.4.1.B' in [Gromov1986].    So essentially  we  re-expose the  Gromov
proof of  Theorem 3.4.1.B' in our  concrete case, and we try to make   Gromov's argument more transparent.
Now to prove the existence of a $\beta^k_{N}$-regular  isometric immersion
$f_3$ which is $C^0$-close to a given map $f_0$, it  suffices to find a section of this extended differential relation  which lies over $f_0$ (Proposition 3.12).  That is only the essential new ingredient in  our
proof of Proposition 3.7.

\medskip

Now we are going to define our extended differential relation.
Let  us denote also  by $f_0$  a map $M^{m+1} \to (\R^{kN}, \beta^k_{N})$ extending 
a given map $f_0 : M^m \to \R ^{kN}$. We denote by $ F_0$ the  corresponding section of the bundle $M^{m+1}  \times \R^{kN} \to M^{m+1}$, i.e.
$F_0 (v) = (v,  f_0 (v))$.
 Denote by $\Gamma _0 \subset M^{m+1}\times  \R^{kN}$
the graph of $f_0$ (i.e. it is the image of $ F_0$), and let $p^*( g)$ and $p^*( \beta^k_{N})$ be the pull-back of the
forms $g$ and  $\beta^k_{N}$ to $M^{m+1} \times  \R^{kN}$ under the obvious projection. Take  a small
neighborhood  $ Y \supset \Gamma_0$ in $M^{m+1}\times \R^{kN}$. Since $\beta^k_{N}$ and $g$ are exact
forms we get

$$ p^*( \beta^k_{N}) - p^*( g) =  d \hat \beta_N$$
for some smooth $(k-1)$-form $\hat \beta_N$ on $ Y$.
\medskip

Our next  observation  is

\medskip

{\bf 3.8. Lemma.} {\it Suppose  that a map $ F : M^{m+1} \to Y$
corresponds to a $\beta^k_{N}$-regular immersion $ f : M^{m+1}\to \R^{kN}$.
Then $F$ is a $d\hat \beta^k_{N}$-regular immersion.}

\medskip

{\it Proof.}  We need to show that for all $y= F (z)\in Y$, $z\in M^{m+1}$, the  composition
$\rho$ of the maps 
$$ T_y Y \stackrel{I_{p^*( \beta^k_{N})-p^*( g)}}{\to} \Lambda ^{(k-1)} T_y Y \to \Lambda ^{(k-1)} (d F(T_{(z)}(M^{m+1}) ) $$
is onto.  This  follows from the consideration of the restriction of
$\rho$ to the subspace $S \subset T_y Y$  which is tangent to the fiber $ \R^{kN}$ in $M^{m+1}\times  \R^{kN} \supset Y$.\QED

\medskip

Now for a map $d\hat\beta_{N}$-regular map $F : M^{m+1} \to  Y$ and a $(k-2)$-form  $\phi$ on $M^{m+1}$ we set 
$$ \Dd ( F, \phi) : =  F^* (\hat \beta_N) + d\phi .\leqno(3.9) $$

With this notation the map $ f : M^{m+1} \to \R^{kN}$ corresponding to $ F: M^{m+1} \to Y$ satisfies 
$$ f^* (\beta_{N}^k) =  F^* (p^*( \beta_{N}^k)) = g +  F ^* (d \hat \beta_N) = g + d \Dd ( F, \phi),$$
 for any $\phi$. Since the space of $(k-2)$-forms $\phi$ is contractible, it follows that the space of $d\hat \beta_{N}$-regular sections  $ F : M^{m+1} \to Y$ for which
$$f^* (\beta_{N}^k) = g + dg_1\leqno (3.9.1)$$
 for a given $(k-1)$-form $g_1$ has the same homotopy type
as the space of solutions to the equation 
$$   \Dd ( F, \phi) =  g_1.$$

In particular  the equation $ f^*_3(\beta_{N}^k) = g$ reduces to the equation
$\Dd ( F, \phi) = 0$ in so far as the unknown map $ f_3$ is $C^0$-close to $f_0$ (so that its graph lies inside $Y$).

\medskip

We define by $\tilde \Phi_{reg}$ the solution sheaf of the equation (3.9) whose component $F$ is $d\hat\beta_{N}$-regular.

\medskip

{\bf 3.10. Proposition.} {\it The  restriction of the solution sheaf  $\tilde \Phi _{reg}$ to $M^{m}$ satisfies the
$H$-principle. Hence the solution sheaf of $\beta_{N}^k$-regular isometric immersions
$f: (M^{m+1},g)\to (\R^{kN}, \beta_{N}^k)$ such that $F(M^{m+1}) \subset Y$ restricted to $M^m$  also  satisfies the
$H$-principle. }

\medskip

Before proving  this Proposition we shall prove the following 

\medskip

{\bf 3.11. Lemma}. {\it The
differential operator $ \Dd$ is infinitesimal invertible at those
pairs $(F, \phi)$ for which the underlying map $ f$ is a $\beta_N^k$-regular
immersion.}

\medskip

{\it Proof.} The  linearization $L_{(F, \phi)}\Dd$ acts on the space of   couples
$(V, \tilde \phi)$ where $V$ is  a section of $f^* (T_*\R ^{kN})$ (a vector field on $\R^{kN}$ along
 the corresponding map $f$) and $\tilde \phi$ is a $(k-2)$-form on $M^{m+1}$ as follows
$$L_{(F , \phi)}\Dd(V, \tilde \phi) =  \Ll_ V \hat \beta_N + d\tilde \phi = V \rfloor d\hat \beta _N+ d ( V\rfloor \hat \beta_N) + d\tilde  \phi.\leqno (3.11.1)$$
 By Lemma 3.8 the map $ F$ is a $d\hat \beta_N$-regular immersion. Hence the equation for $V$
$$F^*(V \rfloor d\hat \beta_{N} )= \tilde g,\leqno(3.11.2)$$
is solvable for all $(k-1) $-form $\tilde g$ on $M^{m+1}$.  Now  we set: 
$$   \tilde \phi := F^*( V \rfloor \hat \beta_{N})\leqno (3.11.3)$$
 Clearly every couple
$(V, \tilde \phi)$ satisfying (3.11.2) and  (3.11.3)  is a solution of the  equation $L_{(F, \phi)} \Dd (V, \tilde \phi)=
\tilde g$ for any given $(k-1)$-form $\tilde g$. \QED

\medskip

{\it Proof of Proposition 3.10}.  Taking into account Lemma 3.11 and 2.3 (Nash implicit function theorem),  2.5 (Nash implicit function theorem implies the microflexibility) we get the microflexibility of $\tilde
\Phi_{reg}$. Next we use the Gromov observation [Gromov1986, 3.4.1.B'] that $M^m$
is a sharply movable submanifold by acting diffeotopies in $M^{m+1}$ which implies that
the restriction of $\tilde \Phi_{reg}$ to $M^m$ is  flexible. Hence we get the first
statement of Proposition 3.10 immediately. The second statement follows by a remark above relating (3.9) and (3.9.1).\QED

\medskip

{\it  Completion of the proof of Proposition 3.7.}

Suppose we are given a map 
$ f_0 : M^m \to \R^{kN}$.
Since $M^m$ is a deformation retract of $M^{m+1}$ the map $f_0$ extends
to a map $f: M^{m+1} \to \R^{kN}$. 

For each $z \in M^{m}$ we denote by $Mono (( T_z M^{m+1} , g), (T_{ f(z)}\R^{kN}, \beta_{N}^k))$
the set of all monomorphisms $\rho: T_zM^{m+1} \to T_{f(z)}\R^{kN}$ such that 
the restriction of $\beta_{N}^k(f(z))$ to $Df ( T_z M^{m+1})$ is equal to $(Df^{-1}) ^*g$.

\medskip

{\bf 3.12. Proposition.} {\it There exists  a section $s$ of the fibration $Mono((TM^{m+1}, g),   f^* (T\R^{kN}, \beta_{N}^k))$ 
such that $s(z) (T_z M^{m+1})$ is $\beta_{N}^k$-regular subspace for
all $z\in M^{m}$.}

\medskip

{\it Proof of Proposition 3.12.} The proof of Proposition 3.12 consists of 3 steps.  

\medskip

\underline{Step 1.}  We consider $TM^{m+1}$ and $M^{m} \times \R^{kN}$ as  vectors bundles
over the same base $M^{m}$. We  shall
show the existence of a section $s_1\in Mono(TM^{m+1},M^{m}\times \R^{kN}) $ such that  the image  $s_1(TM^{m+1})$ is a $\beta^k_{N}$-regular sub-bundle of dimension $(m+1)$ in $M ^{m}\times \R^{kN}$. 
To save notations  we also
denote by $\beta^k$ the following k-form on $\R^{kN}=\oplus _{j=1} ^N \R ^k_j$ 
$$\beta^k =  \sum_{j=1} ^{N} dx_j^1 \wedge  \cdots \wedge dx^k_j.$$
Here  $(x ^i_j), 1\le i \le k, $ are coordinates in $\R^k _j$ for each  $j =\overline{1, N}$.

\medskip
We put for $l\ge k \ge 3$

$$ \delta (l,k) : = (l-1) + {k-1 \over 2} ( 2 + [ { l-2 \over  k-1 } ] ) ([{ l-2 \over k-1 } -1 ] ) + [ {l-1\over k-1} ] ( 1+ ( (l-1) \mod 
 (k-1))  ).$$
 Here  we set  $ i \mod (k-1) := i - (k-1)\cdot [i/(k-1)]$.

 \medskip

{\bf 3.13. Lemma}. {\it For each given $l \ge k\ge 3$ there there exists a l-dimensional
subspace $V^{l}$ in $\R^{kN}$ such that  $V^{l}$ is $\beta^k$-regular subspace, provided that
$N \ge \delta (l,k)$. }

\medskip

{\it  Proof.}  We shall construct a linear embedding $f^l:V^l \to \R^{kN}$ whose image
satisfies the condition of Lemma 3.13.  We work  in  opposite way, i.e. for each $l$ we shall find
a number $\delta(l,k)$  and an embedding $f: V^l \to \R^{k\delta(l,k)}=\oplus _{j=1}^{\delta(l,k)}\R^k_j$  and  $f$ can be written as
$$f:= f^l = (f_1^l , f_2^l, \cdots , f_{\delta(l,k)}^l),\, f_j^l : V^l \to \R^k _j, j =\overline{1,\delta(l,k)},$$
such that $f$ satisfies Lemma 3.13 with $\delta (l,k) = N$. Clearly  the embedding
$V^l \to V^{k \delta(l,k)} \to \R ^{kN}$ also satisfies the condition of Lemma 3.13 for
all $N \ge \delta(l,k)$.

We can assume that $V^k\subset V^{k+1} \subset \cdots \subset V^l$ is a chain
of subspaces in $V^l$ which is generated by some vector basis $(e_1, \cdots ,
e_l)$ in $V^l$. We denote by $(e_1 ^*, \cdots , e_l^*)$ the dual basis
of $(V^l)^*$. By construction, the restriction of $(e_1^*, \cdots , e_i ^*)$ to
$V^i$ is the dual basis of $( e_1, \cdots, e_i) \in V^i$.  For the sake of simplicity we shall denote
the restriction of any $v_j^*$ to these subspaces also by $v_j ^*$ (if the restriction is not zero). We shall  construct $f_i^l$ inductively   on the dimension $l$ of $V^l$ such that the
following condition holds for all $ k\le i \le l$
$$ <(f_1 ^l)^* (\Lambda ^{k-1} (\R^k_1)),  (f_2^l) ^*(\Lambda ^{k-1} (\R^k_2)),\cdots 
, (f_ {\delta(i,k)}^l) ^*(\Lambda ^{k-1} (\R^{k-1}_{\delta(i,k)}))>_{\otimes \R}= \Lambda ^{k-1} (V^i)^*.\leqno (3.14) $$
The condition (3.14) implies that 
$$I_{\beta^k} (\R^{ k\delta(l,k)}) = \Lambda ^{k-1} (V^i) ^*,$$
so $f^i(V^i)$ is $\beta^k$-regular.
For $l=k$ we can take $f_1^k = Id$, and $\delta (k,k) =1$. Suppose that $(f_1^i, \cdots, f_{\delta(i,k)}^i)$ are  already constructed for  our map 
$$f^i = (f_1^i, \cdots, f_{\delta(i,k)}^i) :  V^i \to \R^k_1 \times  \cdots  \times \R^k _{\delta(i,k)}.$$
We shall construct map $f^{i+1}$ as follows. We  set  for $j \le \delta (i,k)$
$$f^{i+1} _j ( e_p) = f^i_j (e_p) \text { if }  1\le p \le  i, $$
$$f^{i+1} _j   (e_{i+1}) =  0 .$$
To find
$f_j^{i+1} $, $ \delta (i,k) +1 \le j \le\delta(i+1,k)$, so that (3.14) holds for the next induction step  $(i+1)$, it suffices to find   linear maps  $f_{\delta (i,k) +1}^{i+1}, \cdots ,
f_{\delta (i+1,k)}^{i+1}$ with the following property 
$$< (f_{\delta(i)+1}^{i+1}) ^* \Lambda^{k-1} (\R ^k_{\delta (i,k)+1})^*, \cdots , (f_{\delta(i+1,k)}^{i+1})^*\Lambda^{k-1}
(\R^k_{\delta (i+1,k)})^*>_{\otimes \R}\supset  \wedge^{k-1}  (V^i)^*\wedge e_{i+1}^*.\leqno (3.15)$$
We shall proceed as follows. Set $\delta(i+1,k): = \delta (i, k) + [ {i\over (k-1)}] +1$.
  Choose for any $1\le j \le \delta (i+1, k)$  a basic $(w^i_j)$,  $1\le i \le k$,  of the space $\R^k_j$. We let 
$$f_j ^{i+1}( e_{i+1}) = w^1 _j \in \R^k_j, \text { if } j  \ge \delta (i,k) +1,$$
$$ f^{i+1}_{\delta(i,k) +1}  (e_1)  = w ^2_{\delta (i,k) +1}, \:  f^{i+1}_{\delta(i,k)+1}  (e_2) = w ^3 _{\delta(i,k)+1}, \cdots   , f^{i+1}(e_{k-1})=w^k_{\delta(i,k)+1},$$
$$f^{i+1} _{\delta (i,k)+2} (e_k) = w^2_{\delta (i,k) +2}, \: f^{i+1}_{\delta (i) +2} (e_{k+1}) = w ^3_{\delta(i,k) +2}, \cdots , f^{i+1}(e_{2(k-1)})=w^k_{\delta(i,k)+2}, $$
$$ \cdots  $$
$$ \cdots, f^{i+1}_{\delta (i+1,k)} (e_{i} ) = w^{i \mod (k-1)  }_{\delta (i+1,k)} .$$
 It is easy to see that the constructed map  $f^{i+1}$ satisfies (3.15) and hence also (3.14).
 Now   using the identity the $\delta (i+1) - \delta (i) = [ {i\over (k-1)}] +1$ we get 
 $$\delta(l,k)  = 1+ (l-k) + \sum_{i  = k}^{l-1} [ { i \over k-1}]$$
 $$ = (l-1) +  2 (k-1) + 3(k-1)+  \cdots + [{l-2\over k-1}] (k-1) + [ {l-1\over k-1} ] (( (l-1) \mod 
 (k-1)) +1 )$$
 $$ = (l-1) + {k-1 \over 2} ( 2 + [ { l-2 \over  k-1 } ] ) ([{ l-2 \over k-1 } -1 ] ) + [ {l-1\over k-1} ] (1+ ( (l-1) \mod 
 (k-1))  ).$$
\QED


We shall consider $M^{m} \times V^{2m+1}$ as a sub-bundle of $M^{m} \times \R^{\delta (2m+1, k)}$ over $M^{m}$.
Next  we shall find
a section $s_1$ for the step 1 by requiring that $s_1$ is a section of the bundle $Mono (TM^{m+1}, M^8\times V^{v})$ of all  fiber  mono-morphisms from $TM^{m+1} $ to  $M^{m} \times V^{v}$. This section exists, since the fiber $ Mono(T_xM^{m+1}, \R^{v})$ is homotopic equivalent to $SO(2m+1)/SO (m)$  which has all  homotopy groups $\pi_j$ vanishing, if $j\le (m-1)$. This completes
the step 1. 

\medskip

\underline{Step 2}. Once a section $s_1$ in Step 1 is specified we put the following form $g_1$ on
$TM^{m+1}_{| M^m}$:
$$ g_1  = g - s_1 ^* (\beta).$$
In this step we show the existence of a section $s_2$  of
the fibration $Hom ((TM^{m+1}, g_1), (M^{m}\times \R^{k \cdot (m+1) \cdot \binom{m+1}{k}},  \beta^k))$ over $M^m$.
Here we consider $ (M^{m} \times \R^ {k \cdot (m+1) \cdot \binom{m+1}{k}},  \beta^k)$ as a fibration over $M^{m}$
and   equipped with  the  $k$-form $\beta^k$  on the fiber $ \R^{k \cdot (m+1) \cdot \binom{m+1}{k}}$. We do not
require that $s_2$ is a monomorphism. 

 Using the Nash trick [Nash1956]  (see also the proof of Proof of Theorem B.1 below) we can find a finite number of open coverings $U_i^j, \: j =\overline{1,(m+1)},$ of $M^{m}$ which satisfy the following properties:
$$ U^j_i \cap U^j_k = \emptyset, \, \forall j = \overline{1,(m+1)}\text { and } i\not= k,\leqno(3.16)$$ 
and moreover $U^j_i$ is diffeomorphic to an open ball for all $i,j$.
Since $U^j_i$ satisfy the condition (3.16), for a fixed $j$ we can embed the union $\hat U^j =\cup_i U^j_i$ into $\R^ {m}$.  Thus for each $j$ on  the union $\hat U^j$ we have
 local coordinates $x^r_j,\, r =\overline{1,(m+1)},\, j = \overline{ 1,(m+1)}$. Using
partition of unity functions  $f_j (z)$ corresponding to $\hat U^j$ we can write 
$$ g_1(z) = \sum_{j =1} ^{m+1} f_j (z)\cdot \sum_{1\le r_1 < r_2 <  \cdots < r_k\le m+1 }  \mu _{j} ^{r_1r_2 \cdots r_k} (z)\cdot dx_{j}^{r_1} \wedge dx _{j}^{r_2}\wedge \cdots \wedge dx_{j}^{r_k},$$
where the last coordinate $x_j^{m+1}$ corresponds to the direction which is transveral
to $T_z \hat U_j$ in $\hat U_j \times (-1, 1) \subset  M^{m+1}$.
We numerate (i.e. find a function  $\theta$ with values in $\N^+$) on the set $\{ (j, r_1r_2\cdots r_k)\}$ of $N_1 = (m+1)\cdot \binom{m+1}{k} $ elements.  Next we find a section $s_2$ of 
the form
$$ s_2(z) =  (\tilde s_1 (z), \cdots , \tilde s_{N_1} (z)), \: \tilde s_q (z)\in Hom (T_z M^{m+1}, \R^k_ q)$$
such that 
$$ \tilde s_{\theta (j,r_1r_2\cdots r_k)} (z)=f_j (z)\cdot  \mu _j ^{r_1r_2\cdots r_k}(z)\cdot  A _{r_1, r_2,\cdots ,  r_k},$$
$$\text{ where }
A_{r_1, r_2, \cdots , r_k} (\p x_{r_l}\in T_z M^{m+1} ) : = \delta _l^i e_i\in \R^k_q.$$
 Here $(e_1, e_2, \cdots,  e_k)$ is  a vector basis
in $\R^k_q$ for $q = \theta (j, r_1r_2\cdots r_k)$ and $(\p x_{r_l})$ a basic  in $T_z M^{ m+1}  $  defined via  embedding $\hat U_j  \to \R^m$ as above.  Clearly the section $s_2$ satisfies
the condition $s_2^* (\beta ^k(z)) = g_1(z)$ for all $z\in M^{m}$. This completes the second step.
 
\medskip

\underline{ Step 3}.  We put
$$ s = ( s_1, s_2),$$
where $s_1$ is the constructed  section in Step 1 and $s_2$ is  the constructed  section in Step 2. Clearly
$s$ satisfies the condition of Lemma 3.13. \QED

\medskip

 Proposition 3.7   now follows from  Proposition 3.10 and Proposition 3.12.\QED

\medskip

{\bf 3.17. Final remark.} We conjecture that  the isometric embedding
map $\tilde f$ in Theorem 3.6 is unique up to homotopy.  It is the case, if  $\phi^k$ is
a closed stable  form on $M^m$ (i.e.  the  orbit of $GL_x (T_x M^m) (\phi^k)$ is dense  in
the space $\Lambda^k_x  (T^* M)$ for all $x \in M^m$, see [LPV2007] for more information).

\medskip

 {\bf Acknowledgement}. I thank  Misha Gromov   for his illuminating
 discussion   and explanation and the referee for numerous  corrections and helpful
 suggestions. This paper has been written during my stay at IHES, Hokkaido University,
 and MPIMIS in Leipzig  in  2006 under a partial  support  by contract RITA-CT-2004-505493 and
   grant  IAA 100190701.

\bigskip
\begin{center}

{\bf Appendix. } \\
in communication with Kaoru Ono\footnote{Department of  Mathematics,  Hokkaido University,
Sapporo 060-0810, Japan} 
\medskip

{\large  A soft proof of the existence of a universal space.}

\end{center}
For readers' convenience, we present here an elementary proof of the following

{\bf A.1. Theorem.} {\it For any given  positive integers $n, k$ there exists  a smooth
manifold $\Mm$ of dimension $N(n,k)$ and a closed differential k-form $\alpha$ 
on it with the following property. For any closed differential k-form $\om$ on  a smooth manifold $M^n$ such that $[\om] \in H^k (M^n, \Z)$ there is a smooth immersion $f: M^n \to \Mm^{N(n,k)}$ such that $f^* (\alpha) = \om$. }

\medskip

{\it Proof.} As in the proof of Theorem 3.6 we  reduce this problem to the existence of an embedding of $M^n$ to the space $\R^{\overline{N_1}}$ with the constant $k$-form $\beta_{\overline{N_1}}$ such that 
the pull-back of $\beta_{\overline{N_1}}$ is equal to a given exact $k$-form $g$. 
Since $g$ is an exact form there exists a ($k-1$)-form $\phi$ on $M^n$ such that 
$d\phi = \om$.

Next we use the Nash trick of  a construction of an open covering $A_i$ on $M^n$ 
$$ M^n = \cup _{i =0} ^{n}  A_i,\leqno (A.2)$$
such that each $A_i$ is  the union of  disjoint open balls $D_{i,j}$, $j=1, \dots, J(i)$ 
on $M^n$.
(Pick a simplicial decomposition of $M^n$ and construct $A_i$ by the induction on $i$.  
Let $D_{0,j}$ be a small coordinate neighborhood of the $j$-th vertex.  
We may assume that they are mutually disjoint.  
Set $A_0=\cup_{j=1}^{J(0)} D_{0,j}$.  
Suppose that $A_0, \dots, A_i$ are defined.  
Let $D_{i+1,j}$ be a small coordinate neighborhood, which contains  
$S^{i+1}_j \setminus \cup_{\ell=0}^i A_{\ell}$, 
where $S^{i+1}_j$ is the $j$-th $i+1$-dimensional simplex.  
We may assume that they are mutually disjoint.  
Set $A_{i+1}=\cup_{j=1}^{J(i+1)} D_{i+1,j}$.  
Hence we obtain desired open sets $A_0, \dots, A_n$.)   

Let $\{ \rho_i \}$ be the partition of unity on $M$ subordinate to the
covering $\{ A_i \}$. We write $\phi(x) =\sum_{i =0}^{n} \rho_i(x) \cdot \phi$. 
Note that $\omega = d\phi = \sum d\phi_i$.  
Clearly
the form $\phi _i = \rho_i (x) \cdot \phi$  has support on $A_i$.

Let $N_1 = \binom{n}{k-1}$ and
$$ \gamma = \sum_{j=1}^{N_1}  x^1_j dx^2_j \wedge \cdots \wedge dx ^k_j.$$
Note that $j=1, \dots, \binom{n}{k-1}$ are in one-to-one correspondence with 
the sequences $1 \leq i_2 < \dots < i_k \leq n$.  

\medskip

{\bf  A.3. Proposition.} {\it  There is an immersion $f_i : A_i \to (\R^{N_1 k}, \gamma)$ such that $f_i^* (\gamma) = \phi_i.$ In particular, $f_i^* d\gamma = d \phi_i$.} 

\medskip

{\it Proof of Proposition A.3.} Since  $A_i$ is  a union of the disjoint balls
$D_{i,j}$ it suffices to prove the existence of immersion $f_i $ on each ball $D=D_{i,j}$.    
Take some coordinate  $(x_1, \cdots , x_n)$  on the ball  $D$.  
We can write the restriction of the $(k-1)$-form $\phi_i$ to $D$ as $\phi$, where
$$ \phi(x) =\sum _{1\le i_2 < \cdots i_k \le n} \lambda_{i_2\cdots i_k} dx ^{i_2} \wedge \cdots \wedge dx ^{i_k}.$$
We  construct map $f_i$ as follows
$$ f_i (x) = ( \dots, f_{i;i_2\cdots i_k}(x), \dots)_{1\le i_2 < \cdots < i_k \le n}\leqno(A.4) $$
where  for $x = (x_1, x_2, \cdots , x_n)$ we put
$$ f_{i;i_2 \cdots i_k} (x) : D \to \R^k_{i_2\cdots i_k}(x^1, x^2,\cdots ,x^k) ,$$
$$ (x_1, \cdots , x_n)\mapsto (x^1=\lambda _{i_2 \cdots i_k} (x),x^2 = x ^{i_2}, \cdots , x^k =x ^{i_k}).$$
Clearly  we have $f_i^*  (\om) = \phi$.  
It is easy to check that $f_i$ is an immersion on $D$.\QED

\medskip
 We shall use cut-off functions $\chi _i$ with support contained in $A_i \subset M$ 
such that $\chi_i = 1$ on the support of $\rho_i$.  
Then $\widetilde{f}_i=\chi_i \cdot f_i$ can be
extended to the whole $M^n$. 
  
\medskip

Now we construct an immersion
$f: M^n\to \R^{\overline{N_1}} = \R^{N_1k(n+1)}$ by setting 
$$ f(x) = (\tilde f_0, \cdots  , \tilde f_{n}).\leqno (B.5)$$
Clearly $f$ is an immersion such that $f^*\alpha = \omega$. 

\medskip

Finally we note that we can choose $f:M^n \to {\R}^{\overline{N_1}}$ such that its image 
is contained in an arbitrary small neighborhood of the origin.   It suffices to construct
immersion $\tilde f_i$ in (A.5)  such that  the image of $\tilde f_i$ is contained in an arbitrary small neighborhood of  the origin.  Since $\tilde  f_i$ is constructed from
immersion of  ball $D_{ij}$ with help of cut-off function $\chi_i$   such that $|\chi_i| \le 1$
we reduce this problem to construct $f_i$ whose image lies in arbitrary  small neighborhood of origin. We do it by  refining a given  covering $D_{ij}$ of $M^n$ and modifying an immersion $f$  satisfying  the condition of Theorem  A.1.

Choose $R>0$ such that the image of $f$ is contained in the $R$-ball centered at 
the origin $O \in {\R}^{\overline{N_1}}$.  
For a given integer $m$, we pick a refinement $\{V_p \}$ 
of the covering $\{D_{i,j}\}_{i,j}$ such that $f_{i(p)}(V_p)$ is contained in a ball of 
radius $1/m^2$.  (The center of the ball may not be the origin.)  
Here $i(p)$ is chosen so that $V_p \subset A_{i(p)}$, i.e., 
there is $j(p)$ such that $V_p \subset D_{i(p),j(p)}$.  
Applying the Nash trick again to refine $\{V_p\}$ so that 
there is an open covering $\{A'_\ell \}$ of $M$ such that 
each of $A'_i$ is a union of some mutually disjoint family of $V_p$'s.  
On $V_p \subset A'_i$ we modify the construction of the mapping $f_{i;i_2, \dots, i_k}$ as follows.  
Using the translation in $x_2, \dots x_k$-coordinates in each 
${\R}^k_{i_2, \dots, i_k}$, we may assume that 
$$f_{i;i_2, \dots i_k}(V_p) 
\subset [-R,R] \times [-1/m^2, 1/m^2] \times \dots \times [-1/m^2, 1/m^2].$$
Now we consider the mapping 
$$\Phi_m:(x_1, x_2, \dots , x_k) \mapsto (\frac{x_1}{m}, m\cdot x_2, x_3, \dots, x_k).$$
Then we find that $\Phi_m^* dx^1 \wedge \dots \wedge dx^k = dx^1 \wedge \dots \wedge dx^k$ and 
$$\Phi_m \circ f_{i;i_2, \dots , i_k} (V_j) \subset [-R/m,R/m] \times [-1/m,1/m] \times 
[-1/m^2, 1/m^2] \times \dots \times [-1/m^2, 1/m^2].$$

Clearly the  modified map $\tilde f_i$ with  components $\tilde f_{i;i_2, \dots , i_k} (V_j) := \Phi_m \circ f_{i;i_2, \dots , i_k} (V_j)$ (see (A.4) for definition of $f_i$)  together with the new refined partition of
$A_i$ as above  has its image  contained  in an arbitrary  small neighborhood of origin  by taking $m$ arbitrary  large. \QED

\medskip


\medskip
 
\begin{thebibliography}{99999}

\bibitem[D-T1958]{D-T1958}{\sc A. Dold and R. Thom}, {\it Quasifaserungen und 
unendliche symmetrische Produkte}, Annals of Math., vol 67, N 2 (1958), 239-281.
\bibitem[D-P1961]{D-P1961}{\sc A. Dold and D. Puppe}, {\it Homologie nicht-additiver Funktoren. Anwendungen}, Annales  Inst. Fourier, 11 (1961), 201-312.
\bibitem[Gromov1986]{Gromov1986} {\sc M. Gromov}, {\it Partial Differential Relations},
Springer-Verlag 1986, also translated in Russian, (1990), Moscow-Mir.
\bibitem[Gromov2006]{Gromov2006} {\sc M. Gromov}, Privat communication.
\bibitem[Nash1956]{Nash1956} {\sc J. Nash}, {\it The embedding problem for
Riemannian Manifolds}, Ann. Math. 63, N1, (1956), 20-63.
\bibitem[LPV2007]{LPV2007}{\sc H.V. Le, M. Panak, J. Vazura}, {\it Manifolds
admitting stable forms},  Comm. Math.  Univ. Carolinae,  (2007, to appear).
\bibitem[Thom1954]{Thom1954} {\sc R. Thom}, {\it Quelques propri\'et\'es globales
des vari\'et\'es  diff\'erentiables,}  Comm. Math. Helv. 28 (1954), 17-86.
\bibitem[Tischler 1977]{Tischler 1977} {\sc  D. Tischler, } {\it Closed 2-forms and
an embedding theorem for symplectic manifolds}, J.Diff. Geom. 12 (1977), 229-235.
\end{thebibliography}
\end{document}